\documentclass[10pt,reqno]{amsart}

\usepackage{amsmath}
\usepackage{amsthm}
\usepackage{amssymb}
\usepackage{amsfonts,mathrsfs}
\usepackage{stmaryrd}
\usepackage{amsxtra}  
\usepackage{mathtools}   %%% \shortintertext
\usepackage{mathabx}   %%% \ll-symbol
\usepackage{csquotes}   %%% \enquote{}-command
\usepackage{hyperref} 
\usepackage{cleveref}   %%% theorem-reference with \cref{} and \Cref{}

\newtheorem{thm}{Theorem}
\newtheorem*{thm_plain}{Theorem}
\newtheorem{lem}[thm]{Lemma}

\newtheorem{prop}[thm]{Proposition}

\newtheorem{cor}[thm]{Corollary}

\theoremstyle{definition}

\theoremstyle{remark}

\newtheorem*{rmk_plain}{Remark}

\newcommand{\N}{\mathbb{N}}
\newcommand{\Z}{\mathbb{Z}}

\newcommand{\R}{\mathbb{R}}
\newcommand{\C}{\mathbb{C}}
\newcommand{\D}{\mathbb{D}}

\newcommand{\T}{\mathbb{T}}

\newcommand{\coloneqq}{\mathrel{\mathop:}=}

\providecommand{\abs}[1]{\lvert#1\rvert}

\providecommand{\ang}[1]{\langle#1\rangle}

\providecommand{\zbar}{{\bar{z}}}

\DeclareMathOperator{\reg}{Reg}
\DeclareMathOperator{\sng}{Sng}

\DeclareMathOperator{\Ord}{Ord}
\DeclareMathOperator{\supp}{Supp}
\DeclareMathOperator{\IM}{Im}

\newcommand{\eps}{\varepsilon}

\begin{document}
\title[Approximation of polynomial hulls]{Approximation of polynomial hulls by analytic varieties: A counterexample}
\author{Tobias Harz}
\address{Universität Bern, Bern, Switzerland}
\email{tobias.harz@unibe.ch}
\date{\today}
\thanks{The author was partially supported by Schweizerische Nationalfonds Grant 200021-207335.}

\begin{abstract}
We construct a connected, compact set $K \subset \C^2$ with the following property: 
there exist points $p \in \hat{K} \setminus K$ such that there does not exist a sequence $\{A_\nu\}$ of analytic sets 
$A_\nu \subset\subset \C^2$ with boundary satisfying $p \in A_\nu$ for every $\nu \in \N$ and $\lim_{\nu\to\infty} bA_\nu \subset K$. 
For every point in $\hat{K} \setminus K$, we explicitly construct a sequence of Poletsky discs, and we compute the weak limit of the
pushforwards of the Green current under these discs.
\end{abstract}
\maketitle

For every compact set $K \subset \C^n$, let
\begin{align*}
 \hat{K} \coloneqq \{z \in \C^n : \abs{f(z)} \le \sup_K\abs{f} \text{ for every } f \in \C[z_1, \ldots, z_n]\}
\end{align*}
denote the polynomially convex hull of $K$. If $A \subset\subset \C^n$ is an analytic set with boundary $bA \subset K$, 
then the maximum principle shows that $A \subset \hat{K}$. It is well-known from the classical examples of Stolzenberg 
\cite{Stolzenberg63} and Wermer \cite{Wermer82} of polynomially convex hulls with no analytic structure that not every 
point $p$ in the additional hull $\hat{K} \setminus K$ can be explained in this way. However, in those examples, 
for each point $p \in \hat{K} \setminus K$ there does exist a sequence of analytic sets $\{A_\nu\}$ with boundary, 
$p \in A_\nu$, such that the Hausdorff limit $\lim_{\nu\to\infty} bA_\nu$ is contained in $K$. Also, many results 
show that $\hat{K} \setminus K$ carries an analytic structure in a generalized sense, see, e.g., Rossi's local maximum principle, 
and the characterizations of polynomially convex hulls due to Poletsky and Duval--Sibony discussed below. 

Drinovec Dronv\v{s}ek and Forstneri\v{c} \cite{DrinovecDrnovsekForstneric12} showed that for any connected, compact set 
$K \subset \C^n$ that is circular, i.e., invariant with respect to the natural action of the unit circle $\T \subset \C$ on $\C^n$,
every point $p \in \hat{K} \setminus K$ admits a sequence of analytic discs containing $p$ and whose boundaries converge 
to a subset of $K$. They also give an example of a non-connected, compact, circular set $K \subset \C^2$ such that not every 
point in $\hat{K} \setminus K$ admits a sequence of analytic discs as above. Porten \cite{Porten17} adapted this example to show that 
there exist connected, compact, non-circular sets in $\C^2$ with the same property.
 
In this note, we will give an example of the following type, that addresses approximation of polynomially convex hulls 
by arbitrary analytic sets instead of analytic discs.

\begin{thm} \label{T:main}
 There exists a connected, compact set $K \subset \C^2$ with the following property: there exists $p \in \hat{K} \setminus K$ 
 such that there does not exist a sequence $\{A_\nu\}$ of analytic sets $A_\nu \subset\subset \C^2$ with boundary satisfying 
 $p \in A_\nu$ for every $\nu \in \N$ and $\lim_{\nu\to\infty} bA_\nu \subset K$.
\end{thm}
 
The situation is different if we relinquish control on an arbitrary small part of the boundary $bA_\nu$. 
Indeed, Poletsky \cite{Poletsky93} proved the following theorem. 

\begin{thm_plain}[Poletsky]
 Let $K \subset \C^n$ be compact. Then $p \in \hat{K}$ if and only if there exists a sequence of  
 holomorphic maps $f_\nu \colon \bar\D \to \C^n$ such that $\lim f_\nu(0) = p$, and such that the following property is satisfied: 
 for every bounded, pseudoconvex Runge domain containing $K$, every open neighborhood $U$ of $K$ and every $\eps > 0$, 
 there exists $\nu_0 \in \N$ such that for $\nu \ge \nu_0$ one has $f_\nu(\bar\D) \subset \Omega$ and 
 $\sigma(\{\zeta \in \T : f_\nu(\zeta) \notin U\}) < \eps$.
\end{thm_plain}

Another characterization of the polynomial hull is due to Duval--Sibony \cite{DuvalSibony95}.
\begin{thm_plain}[Duval--Sibony]
 Let $K \subset \C^n$ be compact. Then $p \in \hat{K}$ if and only if there exists a positive $(1,1)$-current $T$ 
 with $p \in \supp(T)$ such that $dd^c T = \tilde\sigma - \delta_p$, where $\tilde\sigma$ is a representative Jensen measure 
 for (evaluation at) $p$.
\end{thm_plain}

Let us call a sequence $\{f_\nu\}_{\nu=1}^\infty$ of holomorphic maps $f_\nu \colon \bar\D \to \C^n$ as above Poletsky discs (for $K$)
centered at $p$, and let us call a positive $(1,1)$-current $T$ as above a Duval--Sibony current (for $K$) centered at $p$. 
It was shown by Wold \cite{Wold11} that one can obtain a Duval--Sibony current centered at $p \in \hat{K}$ as the weak limit 
$T \coloneqq \lim_{\nu\to\infty} (f_\nu)_\ast G$ of the pushforwards of the Green current $G$ under a sequence of Poletsky discs 
centered at $p$. Here, the Green current is defined as
\begin{align*}
 \ang{G,\alpha} = -\int_\D \log\abs{\zeta} \cdot \alpha \quad\text{for all } \alpha \in \mathcal{E}^{1,1}(\C)\,.
\end{align*}
In the special case of the compact set $K$ of \Cref{T:main}, for every point $p \in \hat{K} \setminus K$ we will explicitly construct
Poletsky discs $\{f_\nu\}$ for $K$ centered at $p$, and we will compute the limit $T \coloneqq (f_\nu)_\ast G$.

\section{The example}
Let $\D \coloneqq \{\zeta \in \C : \abs{\zeta} < 1\}$ and $\T \coloneqq \{\zeta \in \C : \abs{\zeta} = 1\}$.  Set 
\begin{align*}
 K_1 \coloneqq \{(e^{2\pi i\theta},w) \in \T \times \bar{\D} : \abs{w} = \theta, \theta \in [0,1]\}. 
\end{align*}
Then $K_1$ is connected and compact, and $\bar{\D} \times \{0\} \subset \hat{K}_1$.
Indeed, applying the maximum principle on the discs $\{e^{2\pi i\theta}\} \times \D(0,\theta)$, $\theta \in [0,1]$,
shows that $\T \times \{0\} \subset \hat{K}_1$, and another application of the maximum principle implies that 
$\bar\D \times \{0\} \subset \hat{K}_1$. 

\begin{prop}
 Let $p \in \D \times \{0\}$. There does not exist a sequence of analytic sets $A_\nu \subset\subset \C^2$ with boundary 
 such that $p \in A_\nu$ for every $\nu \in \N$ and $\lim_{\nu \to \infty} bA_\nu \subset K_1$.
\end{prop}
\begin{proof}
Assume, in order to get a contradiction, that a sequence $\{A_\nu\}$ with the listed properties exists
(see, e.g., section 14.3 in \cite{Chirka89} for the precise definition of analytic sets with boundary). 
Without loss of generality each $A_\nu$ is purely $1$-dimensional , 
and we claim that there exist a sequence $\{X_\nu\}$ of Riemann surfaces, 
smoothly bounded open sets $D_\nu \subset\subset X_\nu$, and holomorphic maps $\Phi_\nu \colon X_\nu \to \C^2$ such that 
$p \in \Phi_\nu(D_\nu)$ for all $\nu \in \N$ and $\limsup_{\nu\to\infty} \Phi_\nu(bD_\nu) \subset K_1$. 
Indeed, let $\Omega_\nu \subset \C^2$ be open sets such that $A_\nu$ is an 
analytic subset of $\Omega_\nu$, and let $\rho_\nu \colon \Omega_\nu \to \R$ be a smooth exhaustion function, $\nu \in \N$. 
Then, by local finiteness of $\sng A$ and Sard's theorem, for almost all $t \in \R$ the relative boundary of the set 
$\tilde{D}_\nu^t \coloneqq \{q \in A_\nu : \rho_\nu(q) < t\}$ in $A_\nu$ is a compact smooth manifold contained in $\reg A_\nu$. 
For $\{t_\nu\}_{\nu=1}^\infty \subset \R$ converging to $\infty$ fast enough, the sets $\tilde{D}_\nu \coloneqq \tilde{D}_\nu^{t_\nu}$ 
satisfy $p \in \tilde{D}_\nu$ for all $\nu \in \N$ and $\limsup_{\nu\to\infty} b\tilde{D}_\nu \subset K_1$. For every $\nu \in \N$, 
let $\Phi_\nu \colon X_\nu \to A_\nu$ be the normalization of $A_\nu$, i.e., a finite, proper, holomorphic map from a Riemann surface 
$X_\nu$ to $A_\nu$ that is biholomorphic above $\reg A$ (see, e.g., Proposition 6.2 in \cite{Chirka89}). 
Then $D_\nu \coloneqq \Phi_\nu^{-1}(\tilde{D}_\nu)$ has the desired properties.

Let $(z,w)$ denote the coordinates in $\C^2$, let $p=(z_0,w_0)$, and let $U \subset \C^2$ be an open neighborhood of $K_1$. 
After possibly shrinking $U$, we can assume that there exists a holomorphic branch of $\log z$ on $U-z_0$, 
and that $p \notin \pi_z(U)$, where  $\pi_z$ denotes the projection onto the first variable. For every $\nu \in \N$, 
the boundary $bD_\nu$ consists of finintely many smooth curves $\Gamma_1, \ldots, \Gamma_{N(\nu)}$, 
and since $\limsup_{\nu \to \infty} \Phi_\nu(bD_\nu) \subset K_1 \subset U$, there exists $\nu_0 \in \N$ such that 
for every $\nu \ge \nu_0$ and every $k \in \{1, \ldots, N(\nu)\}$ one has $\Phi_\nu(\Gamma_k) \subset U$. 
Fix $\nu \ge \nu_0$, and define $f_\nu \colon X_\nu \to \C$ as $f_\nu(x) \coloneqq \pi_z(\Phi_\nu(x))-z_0$. 
Since, for every $k \in \{1, \ldots, N(\nu)\}$, the image $\Phi_\nu(\Gamma_k)$ is contained in a domain of definition 
of $\log z = \log \circ \pi_z$, 
it follows that
\begin{align*}
 \int_{\Gamma_k} \frac{df_\nu}{f_\nu} = \int_{\Gamma_k} d(\log f_\nu) = 0\,.
\end{align*}
By the argument principle (see, e.g., Theorem 3.1.9 in \cite{Varolin2011}), this proves that
\begin{align*}
 \sum_{x \in D_\nu} \Ord_x(f_\nu) = \frac{1}{2\pi i} \int_{bD_\nu} \frac{df_\nu}{f_\nu} = 0 \,.
\end{align*}
Hence $f_\nu \neq 0$ on $D_\nu$, i.e., $\Phi_\nu(D_\nu) \cap (\{z_0\} \times \C) = \varnothing$, 
which contradicts the fact that $p \in \Phi_\nu(D_\nu)$ for every $\nu \in \N$. 
\end{proof}

We give another variant of the above example, which we will consider in more detail.
Let $\sigma$ denote normalized arclength measure on $\T$, i.e., $\sigma(\T) = 1$.
Fix a nonempty open set $I \subset \T$ such that $\sigma(\bar{I}) < 1$ and $\sigma(E) = 0$, where $E \coloneqq b_\T I$ 
denotes the relative boundary of $I$ in $\T$, and set
\begin{align*}
 K \coloneqq((\T \setminus I) \times \{0\}) \cup (\bar{I} \times \T)\,.
\end{align*}
In order to provide an example of a connected, compact set, consider also the special case 
$I = I_+ \coloneqq \{\zeta \in \T : \IM(\zeta) > 0\}$, and define
\begin{align*}
 K_2 \coloneqq((\T \setminus I_+) \times \{0\}) \cup (\{-1\} \times \bar\D) \cup (\bar{I}_+ \times \T) \,.
\end{align*}
All of the subsequent statements on $K$ remain true, with the same proofs, for $K_2$.

\begin{lem} \label{T:hull}
 We have $\hat{K} = (\bar\D \times \{0\}) \cup (\bar{I} \times \bar\D)$.
\end{lem}
\begin{proof}
As before, it follows from the maximum principle that $(\bar\D \times \{0\}) \cup (\bar{I} \times \bar\D) \subset \hat{K}$.
For the reverse inclusion, observe first that $\hat{K} \subset\bar\D^2$. Let $p = (z_0,w_0) \in \bar\D^2$ 
such that $z_0 \notin \bar{I}$ and $\abs{w_0}>0$. Since $\C \setminus \bar{I}$ is connected, the set $\bar{I}$ 
is polynomially convex, and we can find $P \in \C[z]$ such that $\abs{P} \le 1$ on $\bar{I}$ and $\abs{P(z_0)} > 1/\abs{w_0}$. 
Then $Q(z,w) \coloneqq P(z)w$ satisfies $\abs{Q} \le 1$ on $K$ and $\abs{Q(p)} > 1$.
\end{proof}

Every point $p = (z,w) \in \hat{K} \setminus K$ that lies in $\bar{I} \times \bar\D$ is contained in an analytic disc 
attached to $K$, namely the disc $\{z\} \times \bar\D$. For the remaining points in $\hat{K} \setminus K$, 
we again have, with the same proof as before, the following result.

\begin{prop}
 Let $p \in \D \times \{0\}$. There does not exist a sequence of analytic sets $A_\nu \subset\subset \C^2$ with boundary 
 such that $p \in A_\nu$ for every $\nu \in \N$ and $\lim_{\nu \to \infty} bA_\nu \subset K$.
\end{prop}

We will explicitly construct Poletsky discs for $K$ centered at any point $p \in \hat{K} \setminus K$. 
If $p = (z_0,w_0) \in \bar{I} \times \bar\D$, then we simply take the single disc $f_\nu(\zeta) \coloneqq (z_0,\zeta)$. 
If $p \in \D \times \{0\}$, then we proceed as follows:  
For every real-valued function $f \in L^1(\T)$, let 
\begin{align*}
 P[f](re^{i\theta}) \coloneqq  \frac{1}{2\pi} \int_{-\pi}^\pi \frac{1-r^2}{1 - 2r\cos(\theta - t) + r^2} f(e^{i\theta}) \,dt
\end{align*}
denote its Poisson integral. Then $P[f]$ is harmonic on $\D$, and extends continuously to all points of $\T$ at which $f$ is continuous. 
Let now $u \coloneqq -P[\chi_{\T \setminus \bar{I}}]$, where $\chi_{\T \setminus \bar{I}}$ denotes the characteristic function of 
$\T \setminus \bar{I}$, let $v$ be a harmonic conjugate of $u$, and set $g \coloneqq e^{u+iv}$. Then $g \colon \D \to \D$ is holomorphic
such that $\abs{g}$ extends continuously to $\bar\D \setminus E$ and such that $\abs{g} = 1$ on $I$ and $\abs{g} < 1$ on 
$\T \setminus \bar{I}$. Note that in the special case $I = I_+$ we have the explicit formula
\begin{align*}
 g(\zeta) = \left(i\frac{1-\zeta}{1+\zeta}\right)^{-\frac{i}{\pi}}\,.
\end{align*}
For every $z_0 \in \D$, let
\begin{align*}
 \varphi_{z_0}(\zeta) \coloneqq \frac{z_0 - \zeta}{1-\zbar_0\zeta}\,.
\end{align*}

\begin{prop}
 Let $p = (z_0,0) \in \D \times \{0\}$. Let $g \colon \D \to \D$ be a holomorphic function such that $\abs{g}$ extends continuously 
 to $\bar\D \setminus E$ and such that $\abs{g} = 1$ on $I$ and $\abs{g} < 1$ on $\T \setminus \bar{I}$. 
 Define $\tilde{f}_\nu(\zeta) \coloneqq (\zeta, g^\nu(\zeta))$. For every sequence $\{r_\nu\}_{\nu=1}^\infty \subset (0,1)$ 
 converging to $1$ fast enough, and for $\tau_\nu(\zeta) \coloneqq r_\nu\zeta$, the holomorphic maps $f_\nu \colon \bar\D \to \C^2$ 
 given by $f_\nu = \tilde{f}_\nu \circ \tau_\nu \circ \varphi_{z_0}$ define a sequence of Poletsky discs for $K$ centered at $p$.
\end{prop}
\begin{proof}
It suffices to consider the case $z_0 = 0$. 
Recall first that every bounded pseudoconvex Runge neighborhood of $K$ contains $\hat{K}$, see, e.g., Theorem 4.3.3 in \cite{Hörmander73}, 
and that $\hat{K}$ has a neighborhood basis of bounded pseudoconvex Runge domains, see, e.g., Theorem 1.3.8 in \cite{Stout07}. 
Let $\Omega \subset \C^2$ be a bounded pseudoconvex Runge domain containing $K$, let $U$ be an open neighborhood of $K$, and let $\eps > 0$.
From the properties of $\abs{g}$ it follows that $\lim_{\nu\to\infty} \abs{g^\nu \circ \tau_\nu} = 0$ locally uniformly on 
$\D \cup (\T \setminus \bar{I})$, and that, for $\{r_\nu\}_{\nu=1}^\infty$ converging to $1$ fast enough, 
we have $\lim_{\nu\to\infty} \abs{g^\nu \circ \tau_\nu} = 1$ locally uniformly on $I$. This shows that, for $\nu \in \N$ large enough, 
the sets $f_\nu(\bar\D)$ will be contained in any given neighborhood of $(\bar\D \times \{0\}) \cup (\bar{I} \times \bar\D) = \hat{K}$, 
and thus in $\Omega$. Moreover, if $B \subset \T$ is an open neighborhood of $E$ such that $\sigma(B) < \eps$, then, 
since $\abs{g^\nu \circ \tau_\nu} \to \chi_I$ uniformly on $\T \setminus B$, we have that for $\nu \in \N$ large enough 
$f_\nu(\T \setminus B) \subset U$. Lastly, $f_\nu(0) = (0, g^\nu(0)) \to (0,0)$ for $\nu \to \infty$.
\end{proof}

For every $p \in \hat{K} \setminus K$, we will also explicitly compute the limit $T = \lim_{\nu \to \infty}(f_\nu)_\ast G$ 
with respect to the Poletsky discs $\{f_\nu\}$ constructed above. We first introduce some notation: For every $z_0 \in \D$, let 
\begin{align*}
 g_\D(z_0,\zeta) = -\log \abs{\varphi_{z_0}(\zeta)} = -\log\left|\frac{\zeta-z_0}{1-\zbar_0\zeta}\right|
\end{align*}
denote the Green's function for $\D$, and let 
\begin{align*}
 \omega_\D(z_0,\zeta) = \frac{1-\abs{z_0}^2}{\abs{\zeta-z_0}^2}\,d\sigma(\zeta)
\end{align*}
be the harmonic measure for $\D$. The Green current with respect to $z_0$ is defined as
\begin{align*}
 \ang{G_{z_0},\alpha} = \int_\D g_\D(z_0, \,\cdot\,) \alpha \quad\text{for all } \alpha \in \mathcal{E}^{1,1}(\bar\D)\,.
\end{align*}
For $z_0 = 0$, we simply write $G = G_0$. With the normalization $d^c = \frac{i}{2\pi}(\bar\partial-\partial)$, 
the Green-Riesz representation formula states that
\begin{align*}
 \ang{dd^cG_{z_0},u} = \int_\T u(\zeta) \,d\omega_\D(z_0,\zeta) - u(z_0) \quad\text{for all } u \in \mathcal{E}(\bar\D)\,,
\end{align*}
i.e., $dd^cG_{z_0} = \omega_\D(z_0,\,\cdot\,) - \delta_{z_0}$. 

Let not now $p = (z_0,w_0) \in \hat{K} \setminus K$. If $p \in \bar{I} \times \bar\D$ and $f_\nu(\zeta) =(z_0,\zeta)$, 
then we simply have $T=(\jmath_{z_0})_\ast G$, where $\jmath_{z_0} \colon \D \hookrightarrow \{z_0\} \times \D$ denotes the 
natural inclusion. If $p \in \D \times \{0\}$, then we first observe that for $g = e^{u+iv}$, $u = -P[\chi_{\T \setminus \bar{I}}]$, 
as constructed above, not only $\abs{g}$ but also $g$ extends continuously to $\bar\D \setminus E$. 
This is due to the fact that, while the harmonic conjugate $v$ of $P[f]$ might not extend continuously to points of $\T$ 
at which $f$ is continuous, it does extend continuously to any point of $\T$ at which $f$ is continuously differentiable 
(see, e.g., I.E.2 in \cite{Koosis98}). It then suffices to prove the following result.

\begin{prop} \label{T:limitGreen}
 Let $p = (z_0,0) \in \D \times \{0\}$. Let $g \colon \D \to \D$ be a holomorphic function such that $g$ extends continuously 
 to $\bar\D \setminus E$ and such that $\abs{g} = 1$ on $I$ and $\abs{g} < 1$ on $\T \setminus \bar{I}$. 
 Define $\tilde{f}_\nu(\zeta) \coloneqq (\zeta, g^\nu(\zeta))$. For every sequence $\{r_\nu\}_{\nu=1}^\infty \subset (0,1)$ 
 converging to $1$ fast enough, and for $\tau_\nu(\zeta) \coloneqq r_\nu\zeta$, the holomorphic maps $f_\nu \colon \bar\D \to \C^2$ 
 given by $f_\nu = \tilde{f}_\nu \circ \tau_\nu \circ \varphi_{z_0}$ have the property that the sequence $\{(f_\nu)_\ast G\}$ 
 converges weakly to 
 \begin{align*}
  T \coloneqq (\imath_0)_\ast G_{z_0} + \int_I (\jmath_\zeta)_\ast G \,d\omega_\D(z_0,\zeta)\,,
 \end{align*}
 where $\imath_0 \colon \D \hookrightarrow \D \times \{0\}$ and $\jmath_\zeta \colon \D \hookrightarrow \{\zeta\} \times \D$ 
 are the natural inclusions, i.e.,
 \begin{align*}
    \langle T,\alpha \rangle 
  = \int_{\D \times \{0\}} g_\D(z_0,\,\cdot\,) \alpha 
   +\int_I \left(\int_{\{\zeta\} \times \D} g_\D(0,\,\cdot\,) \alpha \right) d\omega_\D(z_0,\zeta), 
  \quad \alpha \in \mathcal{E}^{1,1}(\C^2)\,.
 \end{align*}
\end{prop}

We first need the following auxiliary statement. The result is certainly well-known, but due to the lack of a good reference, 
we include a proof for the convenience of reading.

\begin{lem} \label{T:averaging}
 Let $\mu \ll \sigma$ be a finite Borel measure on $\T$. If $p_\nu(\zeta) \coloneqq \zeta^\nu$, $\nu \in \N$, then 
 $\lim_{\nu \to \infty} (p_\nu)_\ast\mu = \mu(\T)\sigma$ weakly.
\end{lem}
\begin{proof}
Since $\mu \ll \sigma$, there exists $a \in L^1(\T,\sigma)$ such that $\mu = a\,d\sigma$, and since trigonometric polynomials are dense 
in $L^1(\T,\sigma)$, there exists $\{a_n\}_{n\in\Z} \subset \C$ such that $a(\zeta) = \sum_{n=-\infty}^\infty a_n\zeta^n$ 
with convergence in $L^1$-norm. Hence, for every $k \in \Z$,
\begin{align*}
   \int_{\T} \zeta^k\,d((p_\nu)_\ast\mu) 
 = \int_{\T} \zeta^{k\nu}a(\zeta)\,d\sigma 
 = \sum_{n=-\infty}^\infty\int_{\T}a_n\zeta^{n+k\nu}\,d\sigma 
 = a_{-k\nu} \to \left\{\begin{array}{c@{\,,\;}l} a_0 & k = 0 \\ 0 & k \neq 0 \end{array} \right.\,,
\end{align*}
which shows that $\lim_{\nu \to \infty} (p_\nu)_\ast\mu$ and $a_0\sigma$ coincide on the set of all trigonometric polynomials. 
Appealing again to the density of trigonometric polynomials in $L^1(\T,\sigma)$, and observing that 
$\mu(\T) = \sum_{n=-\infty}^\infty\int_\T a_n\zeta^n \,d\sigma = a_0$, the statement follows.
\end{proof}

\begin{proof}[Proof of of \Cref{T:limitGreen}]
Fix $\{r_\nu\}_{\nu=1}^\infty \subset (0,1)$ converging to $1$ so fast that $g^\nu \circ \tau_\nu \to g^\nu$ for $\nu \to \infty$ 
locally uniformly on $\T \setminus E$. Note that this is always possible, since for every $\eps > 0$ and every $\nu \in \N$ 
there exists $\delta_\nu > 0$ such that $\abs{p^\nu-q^\nu} < \eps$ for every $p,q \in \bar\D^2$ satisfying $\abs{p-q} < \delta_\nu$, 
and since, by assumption on $g$, for every compact set $L \subset \T \setminus E$ and every $r_\nu \in (0,1)$ close enough to $1$ 
we have $\abs{g \circ \tau_\nu - g} < \delta_\nu$ on $L$. 

Let $T_\nu \coloneqq (f_\nu)_\ast G$, $\nu \in \N$, and observe that it suffices to show that 
\begin{align*}
 \lim_{\nu\to\infty}\ang{T_\nu, dd^c u} = \ang{T, dd^c u} \quad\text{for all } u \in \mathcal{E}(\C^2)\,.
\end{align*}
Fix $u \in \mathcal{E}(\C^2)$, and without loss of generality assume that $u \not\equiv 0$ on $\bar{\D}^2$. Then
\begin{align} \label{E:formulaT}
    \ang{T,dd^cu} 
 &= \ang{(\imath_0)_\ast G_{z_0}, dd^cu} + \int_I \ang{(\jmath_\zeta)_\ast G, dd^c u} \,d\omega_\D(z_0,\zeta) \nonumber\\
 &= \int_\T (u \circ \imath_0)(\zeta) \,d\omega_\D(z_0,\zeta) - (u \circ \imath_0)(z_0) \nonumber \\
 &\quad+ \int_I \left( \int_\T (u \circ \jmath_\zeta) \,d\sigma - (u \circ \jmath_\zeta)(0) \right) d\omega_\D(z_0,\zeta) \nonumber \\
 &= \int_{\T \setminus I} (u \circ \imath_0)(\zeta) \,d\omega_\D(z_0,\zeta) \nonumber \\
 &\quad+ \int_I \left( \int_\T (u \circ \jmath_\zeta) \,d\sigma \right)d\omega_\D(z_0,\zeta) - (u \circ \imath_0)(z_0)\,.
\end{align}
On the other hand, we have
\begin{align} \label{E:formulaTnu}
    \ang{T_\nu,dd^cu}
  &= \int_\T (u \circ f_\nu)(\zeta)\,d\sigma(\zeta) - (u \circ f_\nu)(0) \nonumber\\
  &= \int_\T (u \circ \tilde{f}_\nu \circ \tau_\nu)(\zeta)\,d\omega_\D(z_0,\zeta) - (u \circ \tilde{f}_\nu \circ \tau_\nu)(z_0)\,,
\end{align}
where in the last equation we have used that $(\varphi_{z_0})_\ast \sigma = \omega_\D(z_0,\,\cdot\,)$ 
(note that $((\varphi_{z_0})_\ast \sigma)(M) = \omega_\D(0,\varphi_{z_0}^{-1}(M)) = \omega_\D(z_0,M)$ 
for every Borel set $M \subset \T$, see, e.g., Theorem 4.3.8 in \cite{Ransford95}). Let $\varepsilon > 0$. 
Clearly, for $\nu \in \N$ large enough, we have
\begin{align} \label{E:estimateat0}
 \abs{(u \circ\tilde{f}_\nu \circ \tau_\nu)(z_0) - (u \circ \imath_0)(z_0)} < \eps\,.
\end{align}
Since $\abs{g \circ \tau_\nu} \to \abs{g} < 1$ for $\nu \to \infty$ locally uniformly on $\T \setminus \bar{I}$, 
it follows that $\tilde{f}_\nu \circ \tau_\nu \to \imath_0$ for $\nu \to \infty$ locally uniformly on $\T \setminus \bar{I}$. 
Since $f_\nu(\D) \subset \D^2$ for all $\nu \in \N$, it thus follows from uniform continuity of $u$ on $\bar\D^2$ that for 
$\nu \in \N$ large enough
\begin{align} \label{E:estimateonTminusI}
 \left| 
   \int_{\T \setminus \bar{I}} (u \circ \tilde{f}_\nu \circ \tau_\nu)(\zeta) \,d\omega_\D(z_0,\zeta) 
  -\int_{\T \setminus \bar{I}} (u \circ \imath_0)(\zeta) \,d\omega_\D(z_0,\zeta) 
  \right| < \varepsilon\,.
\end{align}
Fix $\delta > 0$ such that $\abs{u(p)-u(q)} < \eps$ for all $p,q \in \bar\D$ with $\abs{p-q} < \delta$. 
Choose pairwise disjoint, open, connected sets $I_1, \ldots, I_N \subset\subset I$ such that, for every $k = 1, \ldots, N$, 
one has $\abs{\zeta-\zeta'} < \delta/2$ for $\zeta,\zeta' \in I_k$, and such that for $B \coloneqq I \setminus \bigcup_{k=1}^N I_k$ 
we have $\omega_\D(z_0,B) < (\sup_{\bar\D^2}\abs{u})^{-1}\eps$. Note that, in particular, the last condition implies that
\begin{align} \label{E:estimateonB}
 \left| 
   \int_{B} (u \circ \tilde{f}_\nu \circ \tau_\nu)(\zeta) \,d\omega_\D(z_0,\zeta) 
  -\int_{B}\left(\int_\T (u \circ \jmath_{\zeta}) \,d\sigma \right)\,d\omega_\D(z_0,\zeta) 
 \right| < 2\eps\,.
\end{align}
Fix points $\zeta_k \in I_k$, $k = 1, \ldots, N$, that will be further specified later on. By the choice of $\{r_\nu\}_{\nu=1}^\infty$, 
we know that $g^\nu \circ \tau_\nu \to g^\nu$ for $\nu \to \infty$ locally uniformly on $\T \setminus E$. Thus, for $\nu \in \N$ 
large enough, we have 
\begin{align*}
  \abs{(\tilde{f}_\nu \circ \tau_\nu)(\zeta) - \jmath_{\zeta_k}(g^\nu(\zeta))} 
= \abs{(\tau_\nu(\zeta), (g^\nu \circ \tau_\nu)(\zeta)) - (\zeta_k, g^\nu(\zeta))} 
< \delta \quad\text{for all } \zeta \in I_k\,.
\end{align*}
Hence, for every $k = 1, \ldots, N$,
\begin{align*}
 \left| 
   \int_{I_k} (u \circ \tilde{f}_\nu \circ \tau_\nu)(\zeta) \,d\omega_\D(z_0,\zeta) 
  -\int_{I_k} (u \circ \jmath_{\zeta_k})(g^\nu(\zeta)) \,d\omega_\D(z_0,\zeta)
 \right| < \omega_\D(z_0,I_k)\eps\,.
\end{align*}
With $p_\nu(\zeta) \coloneqq \zeta^\nu$ and $\mu \coloneqq g_\ast (\omega_\D(z_0,\,\cdot\,)_{|I_k})$, 
we have $\int_{I_k} (u \circ \jmath_{\zeta_k})(g^\nu(\zeta)) \,d\omega_\D(z_0,\zeta) = \ang{(p_\nu)_\ast\mu, u \circ \jmath_{\zeta_k}}$, 
and since we know from \Cref{T:averaging} that $\lim_{\nu \to \infty} (p_\nu)_\ast\mu = \mu(\T) \sigma  = \omega_\D(z_0,I_k)\sigma$ weakly, 
it follows that, for $\nu \in \N$ large enough,
\begin{align*}
 \left| 
   \int_{I_k} (u \circ \jmath_{\zeta_k})(g^\nu(\zeta)) \,d\omega_\D(z_0,\zeta) 
  -\omega_\D(z_0,I_k)\int_\T (u \circ \jmath_{\zeta_k}) \,d\sigma 
 \right| < \textstyle\frac{1}{N}\eps\,.
\end{align*}
Lastly, since $\omega_\D(z_0,\,\cdot\,)$ is absolutely continuous with respect to $\sigma$, 
by the mean value theorem we can choose the points $\zeta_k \in I_k$ in such a way that, for every $k = 1, \ldots, N$,
\begin{align*}
   \omega_\D(z_0,I_k)\int_\T (u \circ \jmath_{\zeta_k})\,d\sigma 
 = \int_{I_k}\left(\int_\T (u \circ \jmath_{\zeta}) \,d\sigma \right)d\omega_\D(z_0,\zeta)\,.
\end{align*}
Combining the last three displayed formulas, and summing up over $k$, we obtain, with $I' \coloneqq \bigcup_{k=1}^N I_k$,
\begin{align} \label{E:estimateonI'}
 \left| 
   \int_{I'} (u \circ \tilde{f}_\nu \circ \tau_\nu)(\zeta) \,d\omega_\D(z_0,\zeta) 
  -\int_{I'}\left(\int_\T (u \circ \jmath_{\zeta}) \,d\sigma \right)d\omega_\D(z_0,\zeta) 
 \right| < 2\eps\,.
\end{align}
Recalling the formulas \eqref{E:formulaT} and \eqref{E:formulaTnu} for $\ang{T, dd^cu}$ and $\ang{T_\nu, dd^cu}$, 
we see that the estimates \eqref{E:estimateat0}--\eqref{E:estimateonI'} show that 
\begin{align*}
 \abs{\ang{T_\nu,dd^cu} - \ang{T,dd^cu}} < 6\eps
\end{align*}
for $\nu \in \N$ large enough, which proves the claim.
\end{proof}
\begin{rmk_plain}
Note that \eqref{E:formulaT} shows that $dd^cT = \tilde\sigma - \delta_p$ with 
\begin{align*}
   \tilde\sigma 
 = \left(\omega_\D(z_0,\,\cdot\,)_{|\T \setminus I} \times \delta_0\right) + \left(\omega_\D(z_0,\,\cdot\,)_{|I} \times \sigma\right)\,,
\end{align*}
and clearly $\tilde\sigma$ is a representative Jesen measure for $p$.
\end{rmk_plain}

\end{document}